\documentclass[11pt]{article}
\usepackage[a4paper, total={6.5in, 9in}]{geometry}
\usepackage{authblk}
\usepackage{hyperref}
\usepackage{enumitem}
\usepackage{amsmath}
\usepackage{amsthm}
\usepackage{amsfonts}
\usepackage{amssymb}
\usepackage{subfig}
\usepackage{graphicx}
\usepackage{epstopdf}
\usepackage{subfig}
\usepackage{multirow}
\setlist[enumerate]{label={\arabic*.}, ref={\arabic*}}
\newtheorem{theorem}{Theorem}[section]

\newtheorem{proof of lemma}[theorem]{Proof of Lemma}

\theoremstyle{definition}

\newtheorem{example}[theorem]{Example}

\numberwithin{equation}{section}
\author[1]{Shubhangini Gupta}
\author[2]{Prashant Sharma}
\author[3,$\ast$]{Tamal Pramanick}
\affil[1,2,3]{\em \scriptsize {Department of Mathematics, National Institute of Technology Calicut, Calicut, Kerala, India.}}
\affil[1]{\em \scriptsize {email id: shubhangini\_p240187ma@nitc.ac.in}}
\affil[2]{\em \scriptsize {email id: prashant\_m220603ma@nitc.ac.in}}
\affil[3]{\em \scriptsize {email id: tamal@nitc.ac.in, $^\ast$Corresponding author}}

\title{Efficient Numerical Evaluation of Triple Integral Using the Euler's Method and Richardson's Extrapolation}
\usepackage{amsmath,amsthm,amssymb,enumerate}
\usepackage{gensymb}
\usepackage[square,sort&compress,comma,numbers]{natbib}

\date{ }
\begin{document}
\maketitle
\begin{abstract}
In this study, we employ Euler's method and Richardson's extrapolation to solve a triple integral, which is then transformed into a third-order initial value problem. Our objective is to resolve the computational challenges associated with triple integration by transforming it into an initial value problem. Euler's method is the fundamental numerical technique for approximating the solution, thereby establishing a baseline for accuracy. The precision of our computations is subsequently improved by employing Richardson's extrapolation to reduce errors systematically. This approach not only illustrates the adaptability of numerical methods in solving intricate mathematical problems, but it also emphasizes the significance of strategic error reduction techniques in enhancing computational outcomes. We present the efficacy of this method in solving triple integrals in an efficient manner through experimentation and analysis, thereby making a significant contribution to the fields of numerical computation and mathematical modeling.\\[0.2cm]
{\bf Keywords:} Triple integral, Euler’s method, Richardson’s extrapolation, Leibnitz method.\\[0.2cm]
{\bf2010 AMS classification:} 65D30, 65D32, 65L06 
\end{abstract}
\section{Introduction}
The calculation of triple integrals presents a fundamental difficulty across numerous scientific and engineering fields, especially during situations with three-dimensional spatial domains. Existing numerical integration approaches, such as quadrature methods, can occasionally be insufficient for dealing with the computing demands and mistakes that come with generating extremely complicated multi-dimensional integrals, as illustrated in references \cite{filho, haber}. It is easy to see that in recent days the complex mathematical models which appears in various domains such as physics, engineering, and applied sciences, frequently require the utilization of multidimensional integrals. Henceforth these models need to encourage the creation of reliable and efficient numerical techniques for multidimensional mathematical integral, as demonstrated by references \cite{James}. 
Triple integrals which include so many complications, might end up with numerical instability. These require the precise and efficient analysis of their evaluation and hence leading to a significant field of study \cite{Haselgrove, Evans}.

Quadrature methods, known as classical approaches to numerical integration, have been widely used for calculating integrals in various dimensions. Quadrature methods are efficient and highly accurate for one-dimensional integrals, such as Gaussian quadrature, Trapezoidal rule or Simpson's rule. However, when extended to multiple dimensions, quadrature methods often require a tensor product approach which leads to an exponential increase in the number of evaluation points and hence computational complexities.
However, with the continued development of mathematical models, it has become apparent that newer methods need to be developed to handle the computational complexities of integrals in higher dimensions. Unfortunately, such approaches usually show inefficient when compared to triple integrals, cf. \cite{stroud}.
A novel approach has been made by Shapeev and Vorozhtsov, Wong  in the field of numerical integration. However, the challenge of triple integrals remains insufficiently studied.

Over the past few years, there has been a change in attention to applying more specialized numerical methods for integral evaluation with several dimensions. Whereas there have been significant strides in the research on the study of double integrals available. The numerical prediction of triple integrals has attracted little consideration, with only a few works targeting this problem. Youngberg's study of alternative approaches to Romberg integration in Rabiei and Saleeby's recent attempts to define curves for triple integration in \cite{Rabiei} are remarkable examples of the restricted amount of studies to apply numerical methods to three-dimensional integrals.
The present research introduces an innovative approach for making a triple integral into a third-order initial value problem, solving some shortcomings in the existing literature. The aim of this project is to deal with the mathematical tasks linked to triple integration via Euler's method, an existing numerical approach when determining solutions to ordinary differential equations \cite{fox}, and refining its accuracy by using the Richardson's extrapolation, see e.g. \cite{lima, Zlatev, prentice}.
Our study is based on the need to utilize error reduction techniques in numerical computations, particularly when dealing with challenging integral problems. As per the authors' knowledge, the study of evaluating the triple integral using Richardson's extrapolation technique is being introduced for the first time in the literature.
Richardson's extrapolation is widely known for its efficiency in reducing errors across all numerical cases, notably problems connected with boundary values.
By modifying the integration limits to avoid singularities, breaking down the integral into uniform separate components, and managing infinite limits \cite{Guichard}, we may overcome the usual difficulties of triple integration, such as numerical instability and discontinuities. By turning the triple integral into a differential equation framework, we offer ways to apply modern numerical methods. Our method offers major possibilities for improving the accuracy and success rate of multidimensional integral evaluations.

Let $f: \mathbb{R}^3\rightarrow \mathbb{R}$ be a real-valued function which is Riemann integrable in all three variable. For convenience, we assume that $f$ is suitably smooth so that all relevant derivatives exist for our analysis.
Assume $W(x)$ is a function defined through a triple integral operation. This mathematical construct allows us to investigate the cumulative effect of integrating a three-dimensional function \( f(x, y, z) \) over a variable domain, where the limits of integration may vary as a function of $x$ . 
We define \( W(x) \) as follows:
\begin{align}
W(x) = \int_{x_0}^{x} \int_{y_0(x)}^{y_1(x)} \int_{z_0(x,y)}^{z_1(x,y)} f(x, y, z) \, dz \, dy \, dx.\label{eq:g1}
    \end{align}
The limits of each inner integral denoted by \( x_0 \), \( x \), \( y_0(x) \), \( y_1(x) \), \( z_0(x,y) \), and \( z_1(x,y) \), are potentially influenced by the variable \( x \). 
Let $M$ denote the exact numerical value of some mathematical object and $K(h)$ be the approximation of $M$ which depends on the adjustable parameter $h$. Hence we can assume the following
\begin{align}
K(h)=M+a_1(x)h+a_2(x)h^{2}+a_3(x)h^{3}+a_4(x)h^{4}+a_5(x)h^{5}+\cdots,\label{eq:g8}
\end{align}
here $a_i$'s are unknown. 
We can approximate $K(h)$ by using the Richardson extrapolation which yields higher order approximations compared to the original approximations. In the original approximations we will use Euler's method, which is also known as a first order Runge-Kutta method to solve an initial value problem.
The Euler’s method is a straightforward and widely used technique for solving ordinary differential equations, and provides the basis for our initial numerical approximation. Despite its simplicity, Euler’s method provides a useful starting point for iterative refinement. However, to enhance the precision of our computations, we employ the Richardson's extrapolation, a powerful technique that reduces the error associated with numerical approximations. 
Richardson’s extrapolation systematically improves the accuracy of numerical solutions by combining estimates at different step sizes effectively reducing the discretization error. We can refine the results obtained from the Euler's method, effectively increasing the order of accuracy and mitigating the effects of truncation errors. This dual approach not only exhibits the versatility of numerical methods in addressing complex mathematical problems but also highlights the critical role of strategic error reduction techniques in improving computational results. This work is arranged as follows. In the next section \ref{eq:g0}, we have discussed how the triple integral can be converted to the differential equation with help of Leibnitz rule, fundamental theorem of calculus and the chain rule, which produced a third order system of differential equations with initial conditions. Section \ref{eq:g2} presents how the Richardson extrapolation is used to assess error in our numerical technique. The corresponding numerical example is provided in this section. Subsequently in section \ref{eq:g3} we have provided the numerical experiment to support our theoretical findings. We have considered a model problem, which assessed error using various step sizes and solved the problem using Euler's and Richardson extrapolation approach. We observe that both values almost coincide with the decreasing step size. We wrap our work with a prospective outlook on the implications of this work in the last part, section \ref{eq:g4}.

\section{Triple Integral Differentiation Using Leibniz and Chain Rules}\label{eq:g0}
This section explains the fundamental theory required to evaluate the triple integral by applying the Leibnitz rule to differentiate under the integral sign. The following is the initial expression for $W(x)$, is given by
\begin{equation}
W(x) = \int_{x_0}^{x} \int_{y_0(x)}^{y_1(x)} \int_{z_0(x,y)}^{z_1(x,y)} f(x, y, z) \, dz \, dy \, dx.\label{eq:g}
\end{equation}
To begin with, we differentiate \eqref{eq:g} with respect to $x$, with an use of the fundamental theorem of calculus. The first derivative of $W(x)$ with respect to $x$ is denoted by $W'(x)$ and defined as
$$ W'(x) = \frac{d}{dx} \int_{x_0}^{x} \left( \int_{y_0(x)}^{y_1(x)} \left( \int_{z_0(x,y)}^{z_1(x,y)} f(x,y,z) \, dz \right) \, dy \right) \, dx, $$
which imply
\begin{align}
W'(x) = \int_{y_0(x)}^{y_1(x)} \left( \int_{z_0(x,y)}^{z_1(x,y)} f(x,y,z) \, dz \right) \, dy.\label{eq:g5}
\end{align}
To compute $W''(x)$, the second derivative of $W(x)$ with respect to $x$, we need to differentiate $W'(x)$ again with respect to $x$. To do this, we apply the Leibnitz rule for differentiation under the integral sign. Using this, we now proceed to compute $W''(x)$.
Assume that the inner integral in (\ref{eq:g5}) can be expressed as:
\begin{align*}
\int_{z_0(x,y)}^{z_1(x,y)} f(x,y,z) \, dz = F(x,y),
\end{align*}
substituting this result into the expression for $W'(x)$ in (\ref{eq:g5}), we have  
\begin{align*}
W'(x)= \int_{y_0(x)}^{y_1(x)} F(x,y) \, dy.
\end{align*}
Now we differentiate $W'(x)$ with respect to $x$:
\begin{align*}
W''(x) &= \frac{d}{dx} \left( \int_{y_0(x)}^{y_1(x)} F(x, y) \, dy \right),
\end{align*}
using the Leibnitz rule, we obtain
\begin{align}
W''(x) &= \int_{y_0(x)}^{y_1(x)} \frac{\partial F}{\partial x}(x,y) \, dy + F(x,y_1(x)) \frac{dy_1(x)}{dx} - F(x,y_0(x)) \frac{dy_0(x)}{dx}.\label{eq:g6}
\end{align}
To find $W'''(x)$, we need to differentiate $W''(x)$ with respect to $x$. 
Now, we differentiate (\ref{eq:g6}) with respect to $x$ to find the expression of $W'''(x)$:
\begin{align*}
W'''(x)= \frac{d}{dx} \left(\int_{y_0(x)}^{y_1(x)} \frac{\partial F}{\partial x}(x,y) \, dy \right)+ \frac{d}{dx} \left( F(x,y_1(x)) \frac{dy_1(x)}{dx} \right) - \frac{d}{dx} \left( F(x,y_0(x)) \frac{dy_0(x)}{dx} \right),
\end{align*}
with an use of chain rule, this yields
\begin{align*}
W'''(x) &= \frac{d}{dx} \left[ \int_{y_0(x)}^{y_1(x)} \frac{\partial F}{\partial x}(x,y) dy \right] + \frac{d}{dx} \left[ F(x,y_1(x)) \right] \frac{dy_1(x)}{dx} + F(x,y_1(x)) \frac{d^2 y_1(x)}{dx^2} \\
&\hspace*{2cm} - \frac{d}{dx} \left[ F(x,y_0(x)) \right] \frac{dy_0(x)}{dx} - F(x,y_0(x)) \frac{d^2 y_0(x)}{dx^2} ,
\end{align*}
Now assume the following set of equations:
\begin{align}
\begin{aligned}
W' &= P, \label{eq:p}\\
W'' &= Q, \\
W'''&= R,
\end{aligned}
\end{align}
with the following set of initial values:
\begin{equation}
\begin{aligned}
 W(x_0)&= \displaystyle \int_{x_0}^{x_0} \int_{y_0(x)}^{y_1(x)} \int_{z_0(x,y)}^{z_1(x,y)} f(x, y, z) \, dz \, dy \, dx= 0,\\ 
W'(x_0)&=P(x_0)=\displaystyle \int_{y_0(x_0)}^{y_1(x_0)} \left( \int_{z_0(x,y)}^{z_1(x,y)} f(x,y,z) \, dz \right) \, dy,\\ 
\text{similarly,}\quad W''(x_0)&= Q(x_0).
\end{aligned}
\end{equation}
which is obtained from equations \eqref{eq:g} and \eqref{eq:p}. It is clearly observed that
\begin{equation}
\begin{aligned}
P'&= Q,\label{eq:g7}\\
P''&= Q' = R(x),
\end{aligned}
\end{equation}
so the third-order system \eqref{eq:p}, can be expressed with the help of \eqref{eq:g7} as 
\begin{equation}
\begin{aligned}
\begin{pmatrix}
W' \\ P' \\ Q' 
\end{pmatrix}
=
\begin{pmatrix}
  P\\Q\\R(x) \label{eq:g9}
\end{pmatrix}.
\end{aligned}
\end{equation}
Now we can apply Euler's method to the above system which will be of the form
\[
\begin{pmatrix}
  W_{i+1} \\ P_{i+1} \\ Q_{i+1}
\end{pmatrix}
=
\begin{pmatrix}
  W_{i} \\ P_{i} \\ Q_{i}
\end{pmatrix}
+h  
\begin{pmatrix}
  P_{i} \\ Q{_i} \\ R{(x_i)}
\end{pmatrix}
\]
We are considering upper limit $x$ as the last node in the set, $h$ is a step size, which denotes spacing between the nodes $\{x_0,x_1,x_2,\ldots,x_i,\ldots, x\}$. Thus, by numerical solution of this system, we get the values of $W({x_i})$ and $W(x)$.

\section{Evaluating Numerical Approximations via Richardson Extrapolation}\label{eq:g2}
The numerical approximation $K(h)$ is a power series in $h$, so by Richardson extrapolation we can combine linearly to obtain a higher approximation. Hence in view of \eqref{eq:g8}, we have
 \[K_0(x,h)\equiv K(x,h) =K\left( x,\frac{h}{2^{0}}\right)=M(x)+a_1(x)h+a_2(x)h^2+a_3(x)h^3+a_4(x)h^4+\ldots,\]
 \[K_1(x,h)\equiv K(x,\frac{h}{2}) =K\left( x,\frac{h}{2^{1}}\right)=M(x)+a_1(x)\frac{h}{2}+a_2(x)\frac{h^2}{4}+a_3(x)\frac{h^3}{8}+a_4(x)\frac{h^4}{16}h+\ldots, \]
similarly,
$$K_2(x,h)\equiv K\left(x,\frac{h}{2^2}\right) =K\left( x,\frac{h}{4}\right),$$
$$K_3(x,h)\equiv K\left(x,\frac{h}{2^3}\right) =K\left( x,\frac{h}{8}\right),$$
$$K_4(x,h)\equiv K\left(x,\frac{h}{2^4}\right) =K\left( x,\frac{h}{16}\right),$$
$$K_5(x,h)\equiv K\left(x,\frac{h}{2^5}\right) =K\left( x,\frac{h}{32}\right),$$
$$K_6(x,h)\equiv K\left(x,\frac{h}{2^6}\right) =K\left( x,\frac{h}{64}\right).$$
For instance, in order to establish $4^{th}$ order method, we find $d_0,\;d_1,\;d_2,\;d_3$ as follows
\[M_4(x,h)\equiv\sum_{s=0}^{3}d_sK_s (x,h)=M(x)+ O(h^4).\]
The value of $d_0,\;d_1,\;d_2,\;d_3$ can be obtained from the following system
\begin{equation}
\label{eq:soe}
\begin{bmatrix}
1 & 1 & 1 & 1\\\vspace{0.1 cm}
1 & \frac{1}{2} &\frac{1}{4} & \frac{1}{8} \\\vspace{0.1 cm}
1 & \frac{1}{4} &\frac{1}{16} & \frac{1}{64} \\\vspace{0.1 cm}
1 & \frac{1}{8} &\frac{1}{64} & \frac{1}{512} \\
\end{bmatrix}
\begin{bmatrix}
d_0\\\vspace{0.1 cm}
d_1\\\vspace{0.1 cm}
d_2\\\vspace{0.1 cm}
d_3\\
\end{bmatrix}
=
\begin{bmatrix}
1\\\vspace{0.1 cm}
0\\\vspace{0.1 cm}
0\\\vspace{0.1 cm}
0\\
\end{bmatrix},
\end{equation}
by calculating we obtain
$$\begin{bmatrix}
\vspace{0.2 cm} d_0 \\
\vspace{0.2 cm} d_1\\ 
\vspace{0.2 cm} d_2\\ 
d_3\\
\end{bmatrix}
=
\left[\begin{array}{r}
\vspace{0.3 cm} -\frac{1}{21}\\
\vspace{0.2 cm} \frac{2}{3}\\ 
\vspace{0.2 cm} -\frac{8}{3}\\ 
\frac{64}{21}\\
\end{array}\right].
$$
Now suppose we want to construct a $6^{th}$ order method, then we have to find $d_0,\;d_1,\;d_2,\;d_3,\;d_4,\;d_5$ such that
$$M_6(x,h)\equiv\sum_{s=0}^{5}d_sK_s (x,h)=M(x)+ O(h^6).$$
To obtain the values of  $d_0,\;d_1,\;d_2,\;d_3,\;d_4,\;d_5$  solve the following system of equations in the matrix form
\begin{align*}
\label{eq:soe1}
\begin{bmatrix}
1 & 1 & 1 & 1 & 1 & 1\\
\vspace{0.2 cm}
1 & \frac{1}{2} &\frac{1}{4} & \frac{1}{8} & \frac{1}{16} & \frac{1}{32} \\
\vspace{0.2 cm}
1 & \frac{1}{4} &\frac{1}{16} & \frac{1}{64} & \frac{1}{256} & \frac{1}{1024} \\ \vspace{0.2 cm}
1 & \frac{1}{8} &\frac{1}{64} & \frac{1}{512} & \frac{1}{4096} & \frac{1}{32768} \\\vspace{0.2 cm}
1 & \frac{1}{16} &\frac{1}{256} & \frac{1}{4096} & \frac{1}{65536} & \frac{1}{1048576} \\\vspace{0.2 cm}
1 & \frac{1}{32} &\frac{1}{1024} & \frac{1}{32768} & \frac{1}{1048576} & \frac{1}{33554432} \\
\end{bmatrix}
\begin{bmatrix}
d_0\\\vspace{0.2 cm}
d_1\\\vspace{0.2 cm}
d_2\\\vspace{0.2 cm}
d_3\\\vspace{0.2 cm}
d_4\\\vspace{0.2 cm}
d_5\\
\end{bmatrix}
=
\begin{bmatrix}
1\\\vspace{0.2 cm}
0\\\vspace{0.2 cm}
0\\\vspace{0.2 cm}
0\\\vspace{0.2 cm}
0\\\vspace{0.2 cm}
0\\
\end{bmatrix},
\end{align*}
which gives
\begin{align*}
\begin{bmatrix}
\vspace{0.2 cm} d_0\\
\vspace{0.2 cm} d_1\\
\vspace{0.2 cm} d_2\\
\vspace{0.2 cm} d_3\\
\vspace{0.2 cm} d_4\\
d_5
\end{bmatrix}
=
\left[\begin{array}{r}
\vspace{0.2 cm}-\frac{1}{9765}\\
\vspace{0.2 cm} \frac{2}{315}\\
\vspace{0.2 cm}-\frac{8}{63}\\
\vspace{0.2 cm}\frac{64}{63}\\
\vspace{0.2 cm}-\frac{1024}{315}\\
\frac{32768}{9765}
\end{array}\right].
\end{align*}
In general, for an $s^{th}$ order method $M_s(x,h)$ can be expressed as
\[ M_s(x, h) = \sum_{n=0}^{s-1} d_n K_n(x, h) = M(x) + O(h^s),\]
where $d_n$'s are coefficients, and the system is given by
\[
\begin{bmatrix}
A_{1,1} & \cdots &\cdots &\cdots & A_{1,m} \\
\vdots & \ddots &\      &\          & \vdots \\
\vdots &\ & A_{i,j} &  &             \vdots \\
\vdots & \ & \  &\ddots & \vdots\\
A_{m,1} & \cdots &\cdots &\cdots  &A_{m,m}
\end{bmatrix}
\begin{bmatrix}
d_{0}\\
\vdots\\
\vdots\\
\vdots\\
d_{m-1}
\end{bmatrix}
=
\begin{bmatrix}
1\\
0\\
\vdots\\
\vdots\\
0
\end{bmatrix},
\]
where in the matrix $[A_{i,j}]$, the elements are given by $\displaystyle A_{i,j} = \left(\frac{1}{2^{i-1}}\right)^{j-1}$.
Using this, we can find the coefficients for 2nd, 3rd and 5th order methods
\[
\begin{bmatrix}
d_0 \\
d_1
\end{bmatrix}
=
\left[\begin{array}{r}
-1 \\
2
\end{array}\right],\qquad
\begin{bmatrix}
\vspace{0.2 cm}d_0 \\
\vspace{0.2 cm}d_1 \\
d_2
\end{bmatrix}
=
\left[\begin{array}{r}
\vspace{0.2 cm}\frac{1}{3} \\
\vspace{0.2 cm}-2 \\
\frac{8}{3}
\end{array}\right],\qquad
\begin{bmatrix}
\vspace{0.2 cm}d_0 \\
\vspace{0.2 cm}d_1 \\
\vspace{0.2 cm}d_2 \\
\vspace{0.2 cm}d_3\\
d_4
\end{bmatrix}
=
\left[\begin{array}{r}
\vspace{0.2 cm}\frac{1}{315} \\
\vspace{0.2 cm}-\frac{2}{21} \\
\vspace{0.2 cm}\frac{8}{9}\\
\vspace{0.2 cm}-\frac{64}{21}\\
\frac{1024}{315}
\end{array}\right].
\]

\subsection{Error Control}
To achieve a desired level of accuracy in the computation of $W(x)$, an appropriate stepsize $h$ can be chosen. Using a moderately small value of $h$, Euler’s method is applied to solve the equation \eqref{eq:g9} at the nodes $\{x_0,x_1,x_2,\cdots,x_i,\cdots, x\}$, for different step sizes needed to construct 
$M_4(x,h)$ and $M_5(x,h)$ through Richardson extrapolation.\smallskip

The computations are then carried out as follows:
\begin{align*}
M_4(x_i,h)-M_5(x_i,h)=\bar{a}_4(x_i)h^4-O(h^5)\approx \bar{a}_4(x_i)h^4,
\end{align*}
$$\implies\bar{a}_4(x_i) =\frac{M_4(x_i,h)-M_5(x_i,h)}{h^4},$$
which is evaluated at each node $\{x_0 = 1,x_1,x_2,\cdots,x_i,\cdots, x = 5\}$ and yielding the error coefficient $\bar{a}_4$ at each node. It is important to note that, the linear combination is used to construct $M_4(x,h)$, hence $\bar{a}_4$ is not necessarily equal to $a_4$. In fact, we have $a_4=64\,\bar{a}_4$.

The following sequence of calculations is then used to determine a new step size consistent with a tolerance term $\delta$, which is user-defined.
$$H=\left(\frac{\delta}{\max \limits_{x_i} |\bar{a}_4(x_i)|}\right)^\frac{1}{4},\quad n=\left\lceil \frac{x-x_0}{H} \right\rceil\quad\text{and}\quad h=\frac{x-x_0}{n}.$$
Using this new step size, the Euler or Richardson algorithm is repeated, resulting in updated values for $M_4(x_i,h)$ and in particular, specifically $M_4(5,h)\approx W(5)$.

\section{Numerical Experiment}\label{eq:g3}
In the present section, we start with an example as follows.
\begin{example}
We consider the following integral:  
\begin{align}
W(x) = \int_{1}^{5} \int_{0}^{x} \int_{0}^{x+y} \sin(x+y+z) \, dz \, dy \, dx.\label{eq:g10}
\end{align}
From the above integral we have, $f(x, y, z) = \sin(x+y+z)$.
The initial conditions we obtain from \eqref{eq:g10} as
$$ W(1) = 0,\quad W'(1) = P(1) = 1.86017\times 10^{-5},\quad W''(1) = P'(1) = Q(1) = 3.19668\times 10^{-3}. $$
Here the third derivative of $W(x)$ is given by
$$ W'''(x) = Q'(x) = R(x) = 8\sin(4x) - 6\sin(2x)+ \sin(x).$$
Now our aim is to compute $W(x)$. For this purpose we choose various step sizes $h$ to get the desired level of accuracy.
In this example, we assume $h=0.01$ and use Euler's method to obtain the solution of \eqref{eq:g9}.
The nodes $\{x_0,x_1,x_2,\cdots,x_i,\cdots, x\}$ are taken from the limit of the outer integral in \eqref{eq:g10}, where the upper limit is considered as the final node.
Here the starting node $x_0=1$ and ending node $x=5$, with the distance between the nodes are taken at an equidistant interval.\smallskip

Here the exact value of the integral \eqref{eq:g10} at the point $x=5$ is $W(5)=0.193269.$ Now in order to find the approximate solution we apply Euler and Richardson scheme in the integral \eqref{eq:g10}.
In the following table (Table \ref{table:1}) we have described the tolerances and the corresponding stepsizes.

\begin{table}[h!]
\small
\centering
\begin{tabular}{|c|c|c|c|c|c|c|c|}
\hline
$\delta$ & $10^{-16}$ & $10^{-14}$ & $10^{-12}$ & $10^{-10}$ & $10^{-8}$ & $10^{-6}$ & $10^{-4}$ \\
\hline
$h$ & $1.01\times10^{-8}$ & $3.36\times10^{-8}$ & $1.01\times10^{-7}$ & $3.36\times10^{-7}$ & $1.01\times10^{-6}$ & $3.36\times10^{-6}$  & $1.01\times10^{-5}$\\
\hline
\end{tabular}
\caption{\footnotesize Tolerance and corresponding stepsizes}
\label{table:1}
\end{table}

We have calculated the value of $M_4(x_i,h)$ with $h=3.36\times10^{-7}$ and hence we have achieved a tolerance of $\sim\; 10^{-10}$ over the entire interval. It is obvious to see that we have determined the value of $W(5)$ approximately (see Table \ref{table:2}), which is given by
$$W(5)\approx M_4(5,3.36\times10^{-7})= 0.1932547129,$$
which is in the range of the exact value. Table \ref{table:2} gives the integral values for Euler's and Richardson's method for various stepsize $h$. From the table it is very much clear that with the decreasing value of the mesh parameter $h$, both the methods approaching to the exact value. We can easily observe that from initial iteration onwards Richardson scheme is providing the approximate integral value which is very close to the exact solution than the Euler scheme. Henceforth Richardson scheme is faster convergence than the Euler scheme.

\begin{table}[h!]
\centering
\begin{tabular}{|c|c| c |c|c|} 
\hline
$h$ ($\delta=10^{-10}$) & Euler's value & Richardson's value& Difference \\
\hline
$1.152\times10^{-6}$ & $0.16655$ & $0.19322$ & $0.02667$ \\
\hline
$1.139\times10^{-6}$ & $0.16694$ & $0.19322$& $0.02628$ \\
\hline
$1.092\times10^{-6}$ & $0.16804$ & $0.19322$ & $0.02518$ \\
\hline
$9.733\times10^{-7}$ & $0.17062$ & $0.19323$& $0.02261$ \\
\hline
$9.143\times10^{-7}$ & $0.17189$ & $0.19323$& $0.02134$ \\
\hline
$3.365\times10^{-7}$ & $0.18527$ & $0.19325$ & $0.00798$ \\
\hline
\end{tabular}
\caption{\footnotesize Integral values computed by Euler's and Richardson's method with different stepsizes.}
\label{table:2}
\end{table}

The integral curves for Euler's and Richardson's schemes are depicted in the following figures (see Figures \ref{fig:enter-label1}\,-\,\ref{fig:enter-label2}) with various stepsizes. In view of the figures it is clear that, both the integral curves computed for Euler's and Richardson's values are coincides as mesh size $h$
decreases.

\subsection*{4.2 CPU Time Comparison}

To evaluate the computational efficiency of the proposed method, we measured the CPU time taken for the numerical computation of a representative triple integral using both Euler's method and its Richardson extrapolated version. The computations were performed on a system equipped with an Intel Core i5 processor and 8 GB of RAM. The results for various step sizes are presented in Table 3.

\begin{table}[h!]
\centering

\begin{tabular}{|c|c|c|c|c|}
\hline
{Step Size $h$} & \multicolumn{2}{c|}{Euler Method} & \multicolumn{2}{c|}{Richardson Method} \\
\cline{2-5}
& CPU Time (sec.) & Approx. Error & CPU Time (sec.) & Approx. Error \\
\hline
0.01   & 0.024 & 2.31e-3 & 0.039 & 1.47e-4 \\
0.005  & 0.043 & 1.13e-3 & 0.068 & 5.49e-5 \\
0.0025 & 0.077 & 5.62e-4 & 0.091 & 3.11e-6 \\
\hline
\end{tabular}
\caption{Comparison of CPU Time and Approximation Error}
\end{table}

As shown in Table 3, Euler’s method and Richardson extrapolation exhibit comparable computational costs. However, Richardson’s method demonstrates significantly faster convergence and yields a substantial reduction in the approximation error. This trade-off between accuracy and computational effort makes Richardson extrapolation a practical choice for applications where high precision is essential.

\begin{figure}
    \centering
    \includegraphics[scale=0.4]{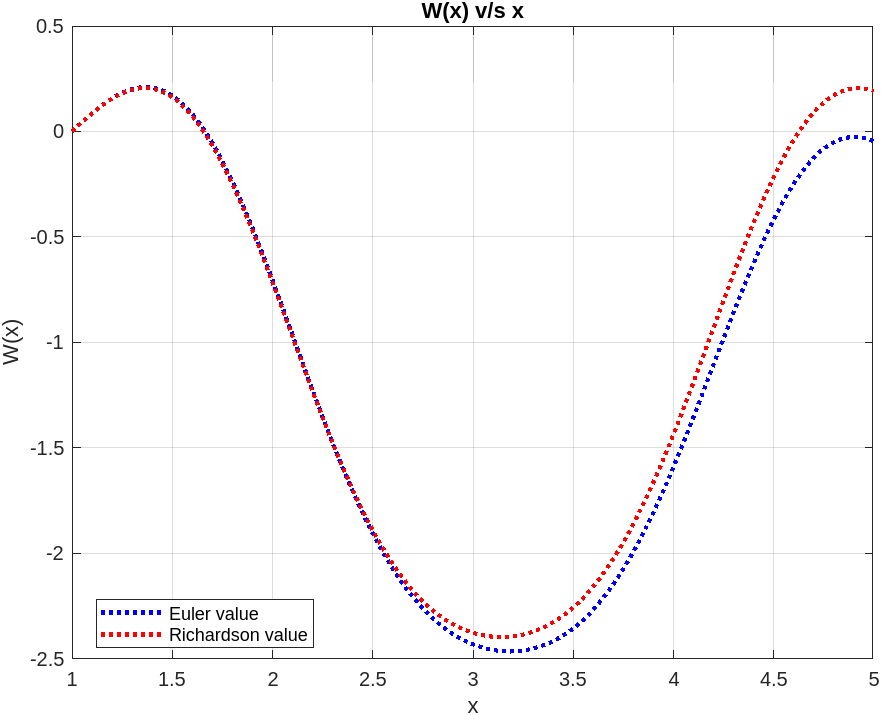}
    \caption{Integral curves for Euler's and Richardson's values at $h=1.13\times10^{-4}$.}
    \label{fig:enter-label1}
\end{figure}
\begin{figure}
    \centering
    \includegraphics[scale=0.4]{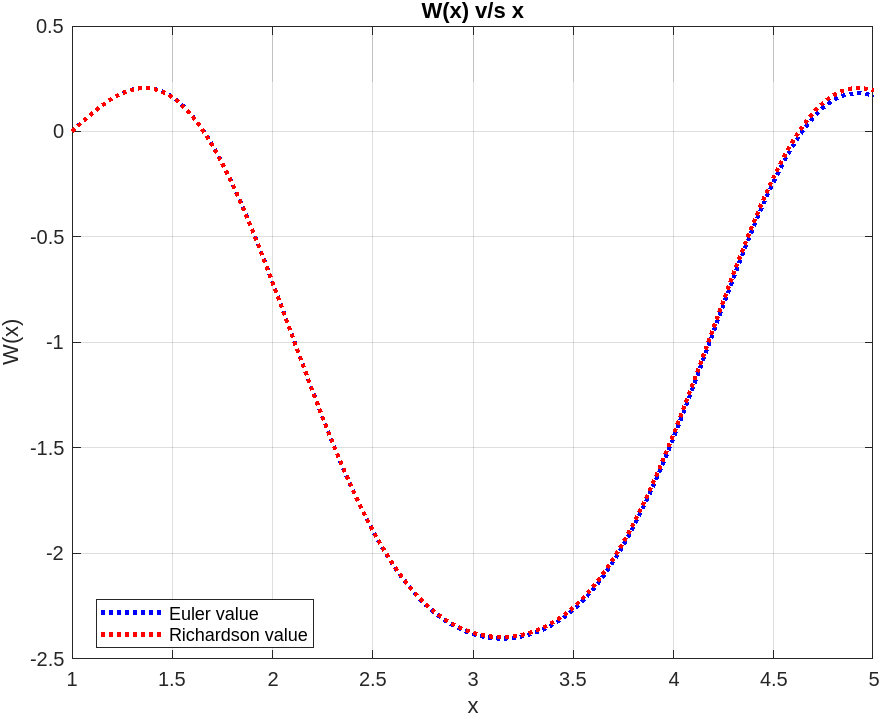}
    \caption{Integral curves for Euler's and Richardson's values at $h=3.36\times10^{-7}$.}
    \label{fig:enter-label2}
\end{figure}

\end{example}
\clearpage

\section{Conclusion}\label{eq:g4}
This study presents a novel method for determining a triple integral and then transitioning to identifying a third-order initial value problem at a significantly later stage.
We achieved exceptional accuracy outcomes by the application of Euler's method and Richardson's extrapolation. However, we achieved an accuracy close to the machine  precision . The methodology employed here is suitable for determining the step-size h necessary for the desired accuracy, which we aim to evaluate based on errors.
If the interval is too big, it's best to break it up into smaller intervals and use the suggested error control approach to get to the goal by combining the intervals in the right way. Also, breaking the interval into smaller intervals and using the process over and over again makes the result more robust and reliable. The suggested method can be extended to higher dimensions, which makes it possible to handle tougher cases like multidimensional integrals. Using this method along with machine learning techniques could also make computation much faster.This approach can be applied to evaluate triple integrals in fluid dynamics, where properties like pressure or velocity fields are integrated over 3D domains. Similarly, it has applications in computational electromagnetics, quantum probability, and even financial models involving three-dimensional parameter spaces.

\section*{Conflict of Interest}

The authors declare that they have no known competing financial interests or personal relationships that could have appeared to influence the work reported in this paper.

\section*{Funding}

This research received no specific grant from any funding agency in the public, commercial, or not-for-profit sectors.

\end{document}